\documentclass[leqno, 11pt, a4paper]{amsart}
\usepackage{amsmath, amsthm, amssymb, mathrsfs, ascmac, bm, enumerate, graphicx, tikz, xcolor, float, thmtools, framed, datetime2}
\everymath{\displaystyle}
\allowdisplaybreaks[4]
\usetikzlibrary{intersections, calc}
\usepackage[geometry]{ifsym}
\DTMlangsetup[en-US]{ord=omit}

\colorlet{headcolour}{black}
\colorlet{rulecolour}{black}

\theoremstyle{definition}

\declaretheorem[title=Theorem, style=plain, numberwithin=section]{thm}

\declaretheorem[title=Corollary, style=plain, sibling=thm]{cor}

\declaretheorem[title=Proposition, style=plain, sibling=thm]{prop}
\declaretheorem[title=Remark, style=definition, sibling=thm]{remark}
\declaretheorem[title=Example, style=definition, sibling=thm]{example}
\declaretheorem[name=Theorem, style=plain]{knownthm}

\numberwithin{equation}{section}

\title[Centralizers of elliptic self-maps]{Centralizers of elliptic univalent self-maps and the embeddability problem in the unit disc}
\author{Akira Mino}
\address{Tohoku University}
\email{mino.akira.t3@dc.tohoku.ac.jp}

\begin{document}

\begin{abstract}
Let $\phi$ be an elliptic univalent self-map of the unit disc in the complex plane that is not an automorphism, with multiplier $0<|\lambda|<1$ at its Denjoy--Wolff point. The simultaneous Schr\"oder linearization of maps commuting with $\phi$ is classical. Using this linearization, we identify its univalent centralizer with the multiplicative semigroup of scalars that preserve the Koenigs domain, and investigate how the geometry of this semigroup controls continuous-time embeddings.

We show that the identity is a non-isolated point of the centralizer precisely when the centralizer contains a non-trivial continuous one-parameter semigroup; such a semigroup need not contain $\phi$. For each logarithmic branch $\mu_k=-\mathrm{Log}\lambda+2\pi i k$, we establish a local sector alternative that characterizes the existence of an embedding with spectral value $\mu_k$.

We prove that the set of admissible branch indices is either empty, a finite interval of integers, or all of $\mathbb{Z}$, and give examples realizing all three possibilities. Finally, we describe the full holomorphic centralizer, including the constant map with value equal to the Denjoy--Wolff point.
\end{abstract}

\maketitle

\section{Introduction}
The \textit{embedding problem} for holomorphic self-maps of the unit disc $\mathbb{D}\subset\mathbb{C}$ asks whether a given map $\phi$ belongs to a continuous one-parameter semigroup $(\phi_t)_{t\geq0}$ at some positive time. Equivalently, after rescaling the time parameter, one may require $\phi_1=\phi$. Since $\phi_{t}\circ\phi=\phi\circ\phi_{t}$ when $\phi_{1}=\phi$, the centralizer of $\phi$ is naturally related to this problem. We write $\mathscr{U}(\mathbb{D})$ for the set of univalent holomorphic self-maps and define the \textit{univalent centralizer} of $\phi$ by
\begin{equation*}
\mathscr{Z}(\phi)=\{\psi\in\mathscr{U}(\mathbb{D}):\psi\circ\phi=\phi\circ\psi\}.
\end{equation*}

Commutativity of elliptic self-maps with a common attracting fixed point has been studied since the work of Pranger [\ref{Pran'70}] and Cowen [\ref{Cow'81},\ref{Cow'84}]. In particular, simultaneous Schr\"oder linearization is classical; see also Cowen's results on the univalence properties of maps commuting with a given elliptic map [\ref{Cow'84}, Corollaries 4.2 and 4.9]. Building on this classical linearization, we use the resulting scalar representation to study the local geometry of the centralizer and the continuous-time embedding problem. For comparison, Contreras, D\'iaz-Madrigal and Gumenyuk [\ref{CDG'26-1}] studied the corresponding problem for non-elliptic univalent maps using additive semigroups in Abel coordinates and obtained related sector alternatives. Section~2 recalls background on non-elliptic centralizers and spirallike functions.

Section~3 develops the multiplicative model and the embedding criteria. The main results of this paper are Theorems \ref{thm:isolation}, \ref{thm:sector} and \ref{thm:admissible-branches}.

Throughout Sections~3 and 4, we consider an elliptic non-automorphism $\phi\in\mathscr{U}(\mathbb{D})$ with Denjoy--Wolff point $\tau\in\mathbb{D}$ and multiplier $\lambda=\phi'(\tau)$ satisfying $0<|\lambda|<1$. Let $h:\mathbb{D}\to\mathbb{C}$ be its normalized Koenigs function, satisfying
\begin{equation*}
h(\tau)=0,\qquad h'(\tau)=1,\qquad h\circ\phi=\lambda h.
\end{equation*}
Write $\Omega=h(\mathbb{D})$ and introduce the multiplicative semigroup
\begin{equation*}
\mathscr{B}_{\phi}^{*}=\{a\in\mathbb{C}^{*}:a\Omega\subset\Omega\}.
\end{equation*}
The classical simultaneous linearization yields the representation $\psi=h^{-1}\circ(a h)$ for $\psi\in\mathscr{Z}(\phi)$ (see [\ref{Pran'70}]). We use the topology of locally uniform convergence on spaces of holomorphic maps and the Euclidean topology on $\mathscr{B}_{\phi}^{*}$. We show that this representation yields a topological semigroup isomorphism $\mathscr{Z}(\phi)\cong\mathscr{B}_{\phi}^{*}$ (Theorem~\ref{thm:model}), which provides a convenient framework for the subsequent results.

Our first result on the topology of the centralizer shows that the identity is a non-isolated point of $\mathscr{Z}(\phi)$ if and only if $\mathscr{Z}(\phi)$ contains a non-trivial continuous one-parameter semigroup (Theorem~\ref{thm:isolation}). However, the semigroup obtained in this way need not have $\phi$ as its time-one map; an example of this phenomenon is given in Example~\ref{ex:not-embedding}.

The main embedding criterion distinguishes the logarithmic branches of the multiplier. Let $\mathrm{Log}$ denote the principal logarithm, and set
\begin{equation*}
\mu_0=-\mathrm{Log}\lambda,\qquad \mu_k=\mu_{0}+2\pi i k,\qquad k\in\mathbb{Z},
\end{equation*}
where $\mathrm{Arg}\lambda\in(-\pi,\pi]$. For each $k$, Theorem~\ref{thm:sector} gives a dichotomy: either an open sector with vertex at $1$ and axis pointing in the direction of $-\mu_k$ is disjoint from
$\mathscr{B}_{\phi}^{*}$, or
\begin{equation*}
 e^{-t\mu_k}\in\mathscr{B}_{\phi}^{*} \qquad\text{for every }t\geq0.
\end{equation*}
In the latter case, the maps $\phi_t^{[k]}=h^{-1}\circ(e^{-t\mu_k}h)$ form a continuous one-parameter semigroup with $\phi_1^{[k]}=\phi$. This provides a local geometric criterion for embeddability along each logarithmic branch, refining the classical criterion formulated in terms of spirallikeness of the Koenigs domain.

We further study the set of admissible branches,
\begin{equation*}
 \mathscr{K}_{\phi}=\{k\in\mathbb{Z}: e^{-t\mu_k}\in\mathscr{B}_{\phi}^{*}\ \text{for every }t\geq0\}.
\end{equation*}
Theorem~\ref{thm:admissible-branches} shows that $\mathscr{K}_{\phi}$ is either empty, a finite interval of integers, or all of $\mathbb{Z}$. Section~4 provides examples realizing all three possibilities, including non-embeddable maps. Finally, in Section~5, we describe the extension from the univalent centralizer to the full holomorphic centralizer. In view of Cowen's classical univalence result, this extension is best regarded as a direct consequence of the common Schr\"oder coordinate rather than as an independent linearization theorem.

\section{Preliminaries}

\subsection{Notation}
Let $\mathbb{D}:=\{z\in\mathbb{C}:|z|<1\}$ denote the unit disc, and set $\mathbb{D}^{*}:=\mathbb{D}\setminus\{0\}$.

 For domains $D,E\subset\mathbb{C}$, we denote by $\mathrm{Hol}(D,E)$ the space of holomorphic maps from $D$ into $E$, and by $\mathscr{U}(D,E)$ the class of univalent maps from $D$ into $E$. When $D=E$, we simply write $\mathrm{Hol}(D)$ and $\mathscr{U}(D)$.

Let $\phi\in\mathscr{U}(\mathbb{D})$ be an elliptic non-automorphism, let $\tau\in\mathbb{D}$ be its Denjoy--Wolff point, and set $\lambda:=\phi'(\tau)$, where $0<|\lambda|<1$. Its normalized Koenigs function is the unique univalent map $h:\mathbb{D}\to\mathbb{C}$ satisfying
\begin{equation*}
h(\tau)=0,\hspace{1em}h'(\tau)=1,\hspace{1em}h\circ\phi=\lambda h.
\end{equation*}
Put $\Omega:=h(\mathbb{D})$. For further background on the Denjoy--Wolff point, Koenigs functions, and continuous one-parameter semigroups, we refer to [\ref{BCD'20}, Chapter 9] and [\ref{CDG'26-1}, Section 2].

\subsection{Non-elliptic self-maps}
For comparison, we recall some results on non-elliptic self-maps due to Contreras, D\'iaz-Madrigal and Gumenyuk.
\begin{knownthm}[{[\ref{CDG'26-1}, Theorem 3.1]}]
Let $\phi\in\mathscr{U}(\mathbb{D})$ be a non-elliptic self-map and let $(S,h,z \mapsto z +1)$ be its canonical holomorphic model. Write $\Omega :=h(\mathbb{D})$ and set $\mathscr{A}_{\phi} :=\{b\in\mathbb{C}: \Omega +b \subset \Omega\}$. Then the following statements hold:
\begin{enumerate}
\item[(A)] $\mathscr{A}_{\phi}$ is a closed additive subsemigroup of $[\mathbb{C},+]$ containing $\mathbb{N}_{0}$.
\item[(B)] If $\phi$ is not an automorphism, then $\mathscr{A}_{\phi} \cap (-\infty,0)=\varnothing$.
\item[(C)] For each $b\in\mathscr{A}_{\phi}$, the map $\psi_{b}(z) := h^{-1}(h(z)+b)$ is well-defined on $\mathbb{D}$ and belongs to $\mathscr{Z}(\phi)$. The map
\begin{equation*}
\mathfrak{T}_{\phi} : \mathscr{A}_{\phi} \to \mathscr{Z}(\phi);\hspace{1em}b \mapsto \mathfrak{T}_{\phi}(b):=\psi_{b}
\end{equation*}
 is injective and continuous, and satisfies $\psi_{b+c} =\psi_{b} \circ \psi_{c}$ for any $b,c \in \mathscr{A}_{\phi}$.
\item[(D)] Moreover, the map $\mathfrak{T}_{\phi}$ is closed, i.e., for every closed subset $X \subset \mathscr{A}_{\phi}$, the image $\mathfrak{T}_{\phi}(X)$ is relatively closed in $\mathrm{Hol}(\mathbb{D})$.
\end{enumerate}
\end{knownthm}

\begin{knownthm}[{[\ref{CDG'26-1}, Proposition 6.9]}]
Let $\phi\in\mathscr{U}(\mathbb{D})$ be a hyperbolic non-automorphism. Then exactly one of the following alternatives holds:
\begin{enumerate}
\item[(i)] $id_{\mathbb{D}}$ is an isolated point of $\mathscr{Z}(\phi)$; in this case, $\phi$ is not embeddable.
\item[(ii)] There exists a unique continuous one-parameter semigroup $(\phi_{t})\subset\mathscr{U}(\mathbb{D})$ with $\phi_{1} =\phi$; in this case, $\mathscr{Z}(\phi)=\{\phi_{t} : t\geq0\}$.
\end{enumerate}
\end{knownthm}

\begin{knownthm}[{[\ref{CDG'26-1}, Theorems 5.4 and 7.9]}]
Let $\phi\in\mathscr{U}(\mathbb{D})$ be a parabolic self-map with Denjoy--Wolff point $\tau\in\partial\mathbb{D}$. Then exactly one of the following alternatives holds:
\begin{enumerate}
\item[(i)] there exist $\rho,\delta>0$ such that
\begin{equation*}
\begin{cases}
\mathscr{A}_{\phi} \cap Sect(\rho,\delta)=\varnothing, & \text{if $\phi$ is of zero hyperbolic step} \\
\mathscr{A}^{*}_{\phi} \cap Sect(\rho,\delta)= \varnothing, & \text{if $\phi$ is of positive hyperbolic step},
\end{cases}
\end{equation*}
where $Sect(\rho,\delta):=\{c\in\mathbb{C}: 0 <|c| <\rho,\ |\arg c| <\delta\}$, $\mathscr{A}_{\phi}:=\{b\in\mathbb{C}: h(\mathbb{D})+b\subset h(\mathbb{D})\}$, $\mathfrak{S}_{\phi}(\psi):=\angle\lim_{z\to\tau}\frac{\psi(z)-z}{\phi(z)-z}\hspace{1em}(\psi\in\mathscr{Z}(\phi))$ and $\mathscr{A}^{*}_{\phi}:=\mathfrak{S}_{\phi}(\mathscr{Z}(\phi))$.

\item[(ii)] there exists a unique continuous one-parameter semigroup
$(\phi_{t})\subset\mathscr{U}(\mathbb{D})$ with $\phi_{1}=\phi$.
\end{enumerate}
\end{knownthm}

\subsection{Spirallike functions}

Let $\mu\in\mathbb C\backslash\{0\}$ with $\mathrm{Re}\,\mu\geq0$. A domain $\Omega\subset\mathbb C$ containing the origin
is called \textit{$\mu$-spirallike} if
\begin{equation*}
e^{-t\mu}\Omega\subset\Omega\qquad\text{for every }t\geq0.
\end{equation*}
A map $h\in\mathrm{Hol}(\mathbb{D})$ is \textit{$\mu$-spirallike with respect to $\zeta\in\mathbb{D}$} if it is univalent, $h(\zeta)=0$ and $h(\mathbb D)$ is $\mu$-spirallike (see [\ref{BCD'20}]). This is the usual notion of \textit{$\mathit{\arg}\mu$-spirallike (with respect to $0$)} in univalent function theory (see [\ref{Pom'75}]).

In particular, when $\mu$ is a positive real number, this condition reduces to the usual starlikeness with respect to $0$.

The following classical characterization provides an analytic criterion for spirallikeness.

\begin{knownthm}[Spirallikeness criterion; see {[\ref{BCD'20}]}]
Let $\tau\in\mathbb{D}$. Let $h\in\mathrm{Hol}(\mathbb{D},\mathbb{C})$ satisfy $h(\tau)=0$ and $h'(\tau)\neq0$, and let $\mu\in\mathbb{C}\backslash\{0\}$ with $\mathrm{Re}\mu\geq0$. Then $h$ is $\mu$-spirallike with respect to $\tau$ if and only if
\begin{equation*}
\mathrm{Re}\left(\frac{(z-\tau)(1-\overline{\tau}z)h'(z)}{\mu h(z)}\right)\geq0\qquad(z\in\mathbb{D}),
\end{equation*}
where the expression at $z=\tau$ is understood by analytic continuation.
\end{knownthm}

Moreover, if $h$ is $\mu$-spirallike with respect to $\tau$, then
\begin{equation*}
\phi_t(z):=h^{-1}(e^{-t\mu}h(z)),\qquad t\geq0,
\end{equation*}
defines a continuous one-parameter semigroup of elliptic self-maps with Denjoy--Wolff point $\tau$, and
\begin{equation*}
\phi_t'(\tau)=e^{-t\mu}.
\end{equation*}
Conversely, the Koenigs domain of every continuous one-parameter semigroup of elliptic self-maps with Denjoy--Wolff point $\tau\in\mathbb{D}$ is $\mu$-spirallike with respect to the origin, where $\mu$ is the spectral value of the semigroup and $\mathrm{Re}\mu\geq0$ (see [\ref{BCD'20}]).

\section{The multiplicative model of the centralizer}

Throughout this section, let
$\phi\in\mathscr{U}(\mathbb{D})$ be an elliptic non-automorphism with Denjoy--Wolff point $\tau$, multiplier $\lambda:=\phi'(\tau)$, normalized Koenigs function $h$, and Koenigs domain $\Omega:=h(\mathbb{D})$. Set
\begin{equation*}
\mathscr{B}_{\phi}^{*}:=\{a\in\mathbb{C}^{*}:a\Omega\subset\Omega\}.
\end{equation*}
We first verify that $\mathscr{B}_{\phi}^{*}$ is closed in $\mathbb{C}^{*}$. If $a_0\notin\mathscr{B}_{\phi}^{*}$, choose $w_0\in\Omega$ with $a_0w_0\notin\Omega$. For all $a$ sufficiently close to $a_0$, the point $a_0w_0/a$ remains in $\Omega$, whereas its product with $a$ does not. Thus the complement of $\mathscr{B}_{\phi}^{*}$ is open.  The Schwarz--Pick lemma applied to $h^{-1}(ah)$ gives $|a|\leq1$.

For later use, we record the basic properties in one statement.
\begin{prop}\label{prop:B}
Let $\phi\in\mathscr{U}(\mathbb{D})$ be an elliptic non-automorphism with Denjoy--Wolff point $\tau$ and multiplier $\lambda:=\phi'(\tau)$. The set $\mathscr{B}_\phi^*$ is a multiplicative subsemigroup of $\mathbb{C}^*$ that is closed in the relative topology on $\mathbb{C}^*$ satisfying $\mathscr{B}_{\phi}^{*}\subset\overline{\mathbb{D}}\setminus\{0\}$.

For every $a,b\in\mathscr{B}_{\phi}^{*}$, the map $\psi_a=h^{-1}\circ(ah)$ belongs to $\mathscr{Z}(\phi)$, and
\begin{equation*}
 \psi_a'(\tau)=a,\hspace{1em}\psi_a\circ\psi_b=\psi_{ab}.
\end{equation*}
Moreover, $\lambda^{-n}\notin\mathscr{B}_{\phi}^{*}$ for every $n\geq1$.
\end{prop}

The proof of this proposition follows the same approach as that of [\ref{CDG'26-1}, Theorem 3.1].


We use the topology of locally uniform convergence on spaces of holomorphic maps and the Euclidean topology on $\mathscr{B}_{\phi}^{*}$.

\begin{thm}\label{thm:model}
Let $\phi\in\mathscr{U}(\mathbb{D})$ be an elliptic non-automorphism. The map
\begin{equation*}
 \mathfrak{S}_\phi:\mathscr{Z}(\phi)\to\mathscr{B}_{\phi}^{*},\hspace{1em}\mathfrak{S}_\phi(\psi)=\psi'(\tau),
\end{equation*}
is an isomorphism of topological semigroups from $[\mathscr{Z}(\phi),\circ]$ onto $[\mathscr{B}_{\phi}^{*},\cdot]$.  Its inverse is
\begin{equation*}
\mathfrak{T}_\phi(a)=h^{-1}\circ(ah),\hspace{1em}a\in\mathscr{B}_{\phi}^{*}.
\end{equation*}
\end{thm}

\begin{proof}
If $\psi\in\mathscr{Z}(\phi)$, then $\psi(\tau)=\tau$ because $\tau$ is the unique fixed point of $\phi$ in $\mathbb{D}$.  Set $a=\psi'(\tau)\neq0$ and $F=(h\circ\psi)/a$.  Then $F\circ\phi=\lambda F$, $F(\tau)=0$, and $F'(\tau)=1$. By the uniqueness of the normalized Koenigs function, we have $F=h$. It follows that $a\Omega=h(\psi(\mathbb{D}))\subset\Omega$, and
$\psi=h^{-1}(ah)$.

Conversely, if $a\in\mathscr{B}_{\phi}^{*}$, then $\psi_a=h^{-1}\circ(ah)$ is a well-defined univalent self-map and
\begin{align*}
h\circ\psi_a\circ\phi&=a\lambda h=\lambda a h \\
&=h\circ\phi\circ\psi_a.
\end{align*}
Hence $\psi_a\in\mathscr{Z}(\phi)$ and $\psi_a'(\tau)=a$.  The map $\mathfrak{S}_\phi$ is continuous because locally uniform convergence of holomorphic functions implies locally uniform convergence of their derivatives.

It remains to prove the continuity of $\mathfrak{T}_{\phi}$. Suppose that $a_n\to a\in\mathscr{B}_{\phi}^{*}$, and let $K\subset\mathbb{D}$ be compact. Since $ah(K)$ is a compact subset of $\Omega$, there exists a compact set $L\subset\Omega$ containing both $ah(K)$ and $a_nh(K)$ for all sufficiently large $n$. The function $h^{-1}$ is uniformly continuous on $L$. Therefore,
\begin{equation*}
h^{-1}(a_n h)\longrightarrow h^{-1}(a h)
\end{equation*}
uniformly on $K$. Hence, $\mathfrak{T}_{\phi}$ is continuous.
\end{proof}

There is also a useful converse construction, which shows that the Koenigs domain is the essential geometric object.

\begin{prop}
Fix $\lambda\in\mathbb{D}\setminus\{0\}$. Let $\mathscr{B}_\phi^*$ is a multiplicative subsemigroup of $\mathbb{C}^*$ that is closed in the relative topology on $\mathbb{C}^*$ such that $1,\lambda\in \mathscr{B}\subset\overline{\mathbb{D}}\setminus\{0\}$. The following are equivalent:
\begin{enumerate}
\item There exists a simply connected domain $\Omega\subset\mathbb{C}$, with $0\in\Omega$, for which $\mathscr{B}=\{a\in\mathbb{C}^{*}:a\Omega\subset\Omega\}$;
\item There exists an elliptic non-automorphism $\phi\in\mathscr{U}(\mathbb{D})$ with $\phi'(\tau)=\lambda$ and $\mathscr{B}_{\phi}^{*}=\mathscr{B}$.
\end{enumerate}
\end{prop}

\begin{proof}
The implication (2) $\Rightarrow$ (1) follows by taking the Koenigs domain of $\phi$. Conversely, let $h:\mathbb{D}\to\Omega$ be a Riemann map and choose $\tau$ with $h(\tau)=0$. Since $\lambda\Omega\subset\Omega$, the formula $\phi=h^{-1}\circ(\lambda h)$ defines an elliptic univalent self-map with fixed point $\tau$ and $\phi'(\tau)=\lambda$. Since $|\lambda|<1$, it is not an automorphism, and the definition gives $\mathscr{B}_{\phi}^{*}=\mathscr{B}$.
\end{proof}

Theorem \ref{thm:model} also yields the following rigidity statement: if $\psi\in\mathscr{Z}(\phi)$ and $0<|\psi'(\tau)|<1$, then $h$ is also the Koenigs function of $\psi$; consequently $\mathscr{Z}(\psi)=\mathscr{Z}(\phi)$.

\begin{thm}\label{thm:isolation}
Let $\phi\in\mathscr{U}(\mathbb{D})\setminus\{id_{\mathbb{D}}\}$ be elliptic. Then $id_{\mathbb{D}}$ is not isolated in $\mathscr{Z}(\phi)$ if and only if $\mathscr{Z}(\phi)$ contains a non-trivial continuous one-parameter semigroup.
\end{thm}

The proof of this theorem follows the same approach as that of [\ref{CDG'26-1}, Theorem 4.4].

\begin{proof}
For an elliptic automorphism, the assertion follows, after conjugation, from the existence of the rotation group. Suppose that $\phi$ is not an automorphism and that $\psi_n\in\mathscr{Z}(\phi)\setminus\{id_{\mathbb{D}}\}$ tends locally uniformly to $id_{\mathbb{D}}$. By Theorem \ref{thm:model}, $h\circ\psi_n=c_nh$, where $c_n\in\mathscr{B}_{\phi}^{*}\setminus\{1\}$ and $c_n\to1$. After taking a subsequence,
\begin{equation*}
 \frac{\log c_n}{|\log c_n|}\to e^{i\theta}.
\end{equation*}
For $t\geq0$, let $m_n=\lfloor t/|\log c_n|\rfloor$. Then $c_n^{m_n}\to e^{te^{i\theta}}$.  Since $\mathscr{B}_{\phi}^{*}$ is a closed multiplicative semigroup, $e^{te^{i\theta}}\in\mathscr{B}_{\phi}^{*}$ for all $t\geq0$.  Therefore,
\begin{equation*}
\Phi_t=h^{-1}\circ\left(e^{te^{i\theta}}h\right),\hspace{1em}(t\geq0)
\end{equation*}
is a non-trivial continuous one-parameter semigroup contained in $\mathscr{Z}(\phi)$. The converse follows immediately from $\Phi_t\to id_{\mathbb{D}}$ as $t\to0^+$.
\end{proof}

Notice that the semigroup obtained in Theorem \ref{thm:isolation} need not contain $\phi$. Embeddability requires a prescribed logarithmic direction, which we now examine.

Let $\phi\in\mathscr{U}(\mathbb{D})$ be an elliptic non-automorphism with Denjoy--Wolff point $\tau$, and let $h$ be its normalized Koenigs function. Suppose that $(\phi_t)_{t\geq0}$ is a continuous one-parameter semigroup satisfying $\phi_1=\phi$. Then there exists $k\in\mathbb Z$ such that
\begin{equation}
h\circ\phi_t=e^{(\mathrm{Log}\phi'(\tau)-2\pi i k)t}h. \notag
\end{equation}
Indeed, by the linearization theorem for elliptic continuous one-parameter semigroups, there exists $\mu\in\mathbb{C}$ such that
\begin{equation*}
h\circ\phi_t=e^{-\mu t}h,\hspace{1em}t\geq0.
\end{equation*}
Since $\phi_1=\phi$ and $h\circ\phi=\phi'(\tau)h$, we have $e^{-\mu}=\phi'(\tau)$. Hence, for some $k\in\mathbb{Z}$, $\mu=-\mathrm{Log}\phi'(\tau)+2\pi i k$.

Let $\mathrm{Log}$ denote the principal logarithm, and define
\begin{equation*}
\mu_0=-\mathrm{Log}\lambda,\qquad\mu_k=\mu_0+2\pi i k,\qquad k\in\mathbb{Z},
\end{equation*}
where $\mathrm{Arg}\lambda\in(-\pi,\pi]$. If a continuous one-parameter semigroup has $\phi$ as its time-one map, then its Koenigs representation must be
\begin{equation*}
 \phi_t^{[k]}=h^{-1}\circ(e^{-t\mu_k}h)
\end{equation*}
for some $k\in\mathbb{Z}$. Thus, the $k$-th branch is admissible precisely when $e^{-t\mu_k}\in\mathscr{B}_{\phi}^{*}$ for all $t\geq0$.

For $v\neq0$ and $\rho,\delta>0$, write
\begin{equation*}
Sect_v(\rho,\delta):=\left\{b\in\mathbb{C}:0<|b-1|<\rho,\left|\arg\frac{b-1}{v}\right|<\delta\right\}.
\end{equation*}

\begin{thm}\label{thm:sector}
Let $\phi\in\mathscr{U}(\mathbb{D})$ be an elliptic non-automorphism. For each $k\in\mathbb{Z}$, exactly one of the following alternatives holds:
\begin{enumerate}
\item there exist $\rho,\delta>0$ such that $\mathscr{B}_{\phi}^{*}\cap Sect_{-\mu_k}(\rho,\delta)=\varnothing$;
\item $e^{-t\mu_k}\in\mathscr{B}_{\phi}^{*}$ for every $t\geq0$, and hence $(\phi_t^{[k]})_{t\geq0}$ is an elliptic continuous one-parameter semigroup with $\phi_1^{[k]}=\phi$.
\end{enumerate}
\end{thm}

The proof of this theorem follows the same approach as that of [\ref{CDG'26-1}, Theorem 5.4].

\begin{proof}
The second alternative excludes the first because $e^{-t\mu_k}=1-t\mu_k+O(t^2)$ as $t\to0^+$.

Suppose the first alternative fails. Choose $a_n\in\mathscr{B}_{\phi}^{*}$ such that
\begin{equation*}
0<|a_n-1|<\frac{1}{n},\hspace{1em}\left|\arg\frac{a_n-1}{-\mu_k}\right|<\frac{1}{n}.
\end{equation*}
Writing $r_n=|a_n-1|$, we have
\begin{equation*}
 \frac{\log a_n}{r_n}\to -\frac{\mu_k}{|\mu_k|},
\end{equation*}
where the logarithm is taken near $1$. Given $t>0$, set $m_n=\lfloor t|\mu_k|/r_n\rfloor$. Then
\begin{equation*}
m_n\log a_n\to -t\mu_k,\hspace{1em}a_n^{m_n}\to e^{-t\mu_k}.
\end{equation*}
Closedness and multiplicative stability of $\mathscr{B}_{\phi}^{*}$ imply $e^{-t\mu_k}\in\mathscr{B}_{\phi}^{*}$. Hence, $\phi_t^{[k]}$ is well defined. The semigroup property and continuity with respect to the topology of locally uniform convergence follow from its Koenigs representation, and $e^{-\mu_k}=\lambda$ gives $\phi_1^{[k]}=\phi$.
\end{proof}

\begin{thm}\label{thm:admissible-branches}
Let $\phi\in\mathscr{U}(\mathbb{D})$ be an elliptic non-automorphism and set
\begin{equation*}
\mathscr{K}_\phi:=\{k\in\mathbb{Z}:e^{-t\mu_k}\in\mathscr{B}_{\phi}^{*}\text{ for all }t\geq0\}.
\end{equation*}
Then $\mathscr{K}_{\phi}$ is either $\varnothing$, $\{m,m+1,\ldots,n\}\ (m\leq n)$ or $\mathbb{Z}$. Moreover, each of these possibilities occurs.
\end{thm}

\begin{proof}
Write $\mu_k=\mu_{0}+2\pi ik$, where $\mathrm{Re}\mu_{0}=-\log|\phi'(\tau)|>0$. If $k_1,k_2\in \mathscr{K}_{\phi}$ and $k_1\le j\le k_2$, write $j=\theta k_1+(1-\theta)k_2$ with $0\le\theta\le1$. Then $\mu_j=\theta\mu_{k_1}+(1-\theta)\mu_{k_2}$, and hence
\begin{equation*}
e^{-t\mu_j}\Omega=e^{-\theta t\mu_{k_1}}e^{-(1-\theta)t\mu_{k_2}}\Omega\subset\Omega.
\end{equation*}
Thus $j\in \mathscr{K}_{\phi}$, so $\mathscr{K}_{\phi}$ is an interval in $\mathbb{Z}$.

Suppose that $\mathscr{K}_{\phi}$ is unbounded above. Choose $k_j\in \mathscr{K}_{\phi}$ such that $k_j\to+\infty$. Then the infinitesimal generators
\begin{equation*}
G_{k_j}(z)=-\mu_{k_j}\frac{h(z)}{h'(z)}
\end{equation*}
satisfy
\begin{equation*}
\frac{G_{k_j}(z)}{2\pi k_j}\longrightarrow -i\frac{h(z)}{h'(z)}
\end{equation*}
locally uniformly in $\mathbb{D}$. By the closedness of the cone of infinitesimal generators, $-ih/h'$ is also a generator. Its associated semigroup satisfies
\begin{equation*}
h\circ\psi_t=e^{-it}h.
\end{equation*}
Thus $\Omega$ is invariant under every rotation, and therefore $\Omega$ is a disc centered at the origin. Consequently, $\mathscr{K}_{\phi}=\mathbb{Z}$.

The case in which $\mathscr{K}_{\phi}$ is unbounded below is analogous. Indeed, choose $k_j\in \mathscr{K}_{\phi}$ such that $k_j\to-\infty$.
Since $-2\pi k_j>0$, the functions $G_{k_j}/(-2\pi k_j)$ are infinitesimal generators and converge locally uniformly to $ih/h'$. Hence $ih/h'$ is also an infinitesimal generator. Its associated semigroup satisfies
\begin{equation*}
h\circ\psi_t=e^{it}h.
\end{equation*}
Thus $\Omega$ is invariant under every rotation and hence is a disc centered at the origin. Consequently, $\mathscr{K}_{\phi}=\mathbb{Z}$.

Therefore, if $\mathscr{K}_{\phi}\neq\mathbb{Z}$, then $\mathscr{K}_{\phi}$ is bounded both above and below. Since it is an interval in $\mathbb{Z}$, it is a finite interval.

The three possibilities are realized by Examples \ref{ex:empty}, \ref{ex:finite-interval} and \ref{ex:branches}, respectively.
\end{proof}

Different admissible branches yield embedding semigroups with different angular velocities in Koenigs coordinates.

Moreover, any two continuous one-parameter semigroups having $\phi$ as their time-one map commute with each other. Indeed, let $h$ be the Koenigs function of $\phi$. Since $\mathscr{Z}(\phi)$ is commutative, for any $t,s\geq0$ and any $k,j\in \mathscr{K}_{\phi}$,
\begin{equation*}
\phi_{t}^{[k]}\circ\phi_{s}^{[j]}=\phi_{s}^{[j]}\circ\phi_{t}^{[k]}.
\end{equation*}

\section{Examples}

\begin{example}\label{ex:not-embedding}
The following example shows that the existence of a non-trivial semigroup in the centralizer does not imply embeddability.  For each $n\in\mathbb{N}\cup\{0\}$, set
\begin{equation*}
\Gamma_n=\{e^n e^{(1+i)s}:s\geq0\},\quad\Gamma=\bigcup_{n=0}^{\infty}\Gamma_n,\quad\Omega=\mathbb{C}\setminus\Gamma.
\end{equation*}
The spirals are pairwise disjoint. Consider the homeomorphism $H:\mathbb{C}\to\mathbb{C}$ defined by $H(0)=0$ and $H(z)=ze^{-i\log|z|}$ for $z\neq0$. It maps each \(\Gamma_n\) onto a radial slit. Since $H(\Omega)$ is starlike with respect to $0$, $\Omega$ is simply connected. Moreover, $e\Gamma_n=\Gamma_{n+1}$ for every $n\in\mathbb{N}\cup\{0\}$, so $e^{-1}\Omega\subset\Omega$. Choose a Riemann map $h:\mathbb{D}\to\Omega$ satisfying $h(0)=0$, and define $\phi=h^{-1}\circ(e^{-1}h)$.  Since $e^{(1+i)t}\Gamma\subset\Gamma$, we have $e^{-(1+i)t}\Omega\subset\Omega$, so
\begin{equation*}
\psi_t=h^{-1}\circ(e^{-(1+i)t}h)\qquad t\geq0.
\end{equation*}
The family $(\psi_{t})_{t\geq0}$ forms a non-trivial continuous one-parameter semigroup in $\mathscr{Z}(\phi)$. Thus, the identity is not isolated in $\mathscr{Z}(\phi)$.

Nevertheless, $\phi$ is not embeddable. The possible spectral values of an embedding semigroup are $\mu_k=1+2\pi i k$, $k\in\mathbb{Z}$. Suppose that the $k$-th branch is admissible. Then $e^{t\mu_k}\Gamma\subset \Gamma$ for all $t\geq0$. Since $1\in \Gamma$, we have $e^{t\mu_k}\in \Gamma$. For $0<t<1$,
\begin{equation*}
|e^{t(1+2\pi ik)}|=e^t<e.
\end{equation*}
Hence $e^{t\mu_k}\in\Gamma_0$. Thus, there exists $s_t\geq0$ such that $e^{t\mu_k}=e^{(1+i)s_t}$. Comparing absolute values gives $s_t=t$, while comparing arguments yields $(2\pi k-1)t\in2\pi\mathbb{Z}$. Since this holds for all sufficiently small $t>0$, we obtain $2\pi k=1$, a contradiction.
\end{example}

\begin{example}\label{ex:branches}
The number of admissible logarithmic branches can vary considerably. If
\begin{equation*}
\Omega=\mathbb{D},\hspace{1em}h=id_{\mathbb{D}},\hspace{1em}\phi(z)=e^{-1}z,
\end{equation*}
then $\mathscr{B}_{\phi}^{*}=\overline{\mathbb{D}}\setminus\{0\}$. Consequently, every $k\in\mathbb{Z}$ is admissible, and
\begin{equation*}
\phi_t^{[k]}(z)=e^{-t(1+2\pi i k)}z
\end{equation*}
defines a continuous one-parameter semigroup satisfying $\phi_1^{[k]}=\phi$. Thus, $\mathscr{K}_{\phi}=\mathbb{Z}$.
\end{example}

The preceding example leads to the following characterization.
\begin{cor}
Let $\phi\in\mathscr{U}(\mathbb{D})$ be an elliptic self-map with Denjoy--Wolff point $\tau\in\mathbb{D}$ and multiplier $\lambda=\phi'(\tau)$, where $0<|\lambda|<1$. Let $h$ be the normalized Koenigs function of $\phi$, satisfying $h(\tau)=0$ and $h'(\tau)=1$. Then $\mathscr{K}_{\phi}=\mathbb{Z}$ if and only if $h(\mathbb{D})=\{z:|z|<1-|\tau|^{2}\}$.
\end{cor}

\begin{proof}
Suppose that $\mathscr{K}_{\phi}=\mathbb{Z}$. By Theorem~\ref{thm:admissible-branches} and its proof, $\Omega=h(\mathbb{D})$ is a disc centered at the origin. Since $h$ maps $\mathbb{D}$ conformally onto a disc centered at the origin and $h(\tau)=0$, there exists $R>0$ such that
\begin{equation*}
h(z)=R\frac{z-\tau}{1-\overline{\tau}z}.
\end{equation*}
The normalization $h'(\tau)=1$ gives $R=1-|\tau|^2$. Hence
\begin{equation*}
h(\mathbb{D})=\{w\in\mathbb{C}:|w|<1-|\tau|^2\}.
\end{equation*}

Conversely, suppose that $h(\mathbb{D})$ is this disc. For every $k\in\mathbb{Z}$ and $t\geq0$, we have $|e^{-t\mu_k}|=|\lambda|^t\leq1$. Thus
\begin{equation*}
e^{-t\mu_k}h(\mathbb{D})\subset h(\mathbb{D})\quad\text{and}\quad\mathscr{K}_{\phi}=\mathbb{Z}.
\end{equation*}
\end{proof}

\begin{example}\label{ex:empty}
The case $\mathscr{K}_\phi=\varnothing$ can occur. Let
\begin{equation*}
h(z)=\frac{z}{1-z^2},\qquad\phi=h^{-1}\circ(\lambda h).
\end{equation*}
The function $h$ is univalent on $\mathbb{D}$. Indeed, if $h(z)=h(w)$ for $z,w\in\mathbb{D}$, then $(z-w)(1+zw)=0$, and hence $z=w$. Moreover, $h(\mathbb{D})=\mathbb{C}\setminus\left(i\left[1/2,\infty\right)\cup(-i)\left[1/2,\infty\right)\right)$. Hence, $h(\mathbb{D})$ is starlike with respect to $0$. Fix a real number $\lambda$ with $0<|\lambda|<1$. Then $\lambda h(\mathbb{D})\subset h(\mathbb{D})$. Therefore, the map $\phi=h^{-1}\circ(\lambda h)$ is a well-defined elliptic univalent self-map of $\mathbb{D}$ with $\phi'(0)=\lambda$. Let
\begin{equation*}
\mu_k=-\mathrm{Log}\lambda+2\pi ik,\qquad k\in\mathbb{Z}.
\end{equation*}
Here, $\mathrm{Log}$ denotes the principal logarithm. For this function $h$, we have
\begin{equation*}
\frac{zh'(z)}{h(z)}=\frac{1+z^{2}}{1-z^{2}}\quad\text{and}\quad\left\{\frac{1+z^{2}}{1-z^{2}}:z\in\mathbb{D}\right\}=\{w\in\mathbb{C}:\mathrm{Re}w>0\}.
\end{equation*}
Since $zh'(z)/h(z)$ maps $\mathbb{D}$ onto the right half-plane, the spirallikeness criterion holds for $\mu_k$ if and only if $\mu_k$ is a positive real number.

If $0<\lambda<1$, then $\mu_k=-\log\lambda+2\pi ik$. Thus, $\mu_k$ is a positive real number if and only if $k=0$. Consequently, $\mathscr{K}_\phi=\{0\}$.

If $-1<\lambda<0$, then $\mathrm{Log}\lambda=\log|\lambda|+\pi i$ and hence $\mu_k=-\mathrm{Log}|\lambda|+(2k-1)\pi i$. For any $k\in\mathbb{Z}$, $\mathrm{Im}\mu_k=(2k-1)\pi\neq0$. For each $k\in\mathbb{Z}$, we can choose $w$ in the right half-plane such that $\mathrm{Re}(w/\mu_k)<0$. Therefore, no logarithmic branch $\mu_k$ satisfies the spirallikeness criterion, and hence $\mathscr{K}_\phi=\varnothing$.
\end{example}

For $m<n$, the domains constructed below are closely related to the standard strongly spirallike domains studied by Sugawa (see [\ref{Sug'12}, p. 27]).

\begin{example}\label{ex:finite-interval}
Every finite interval of integers occurs. Fix integers $m\leq n$, and put $\mu_k:=1+2\pi i k$ for $k\in\mathbb{Z}$. Define
\begin{align*}
&C_{m,n}:=\{s\mu_m+t\mu_n:s,t\ge0\}, \\
&F_{m,n}:=\exp(C_{m,n})=\{e^{x+iy}:x\geq0,2\pi mx\leq y\leq 2\pi nx\}, \\
&\Omega_{m,n}:=\mathbb C\setminus F_{m,n}.
\end{align*}
We first show that $F_{m,n}$ is closed. Let $(w_j)$ be a sequence in $F_{m,n}$ converging to $w\in\mathbb{C}$. Write $w_j=e^{x_j+iy_j}$, where $x_j\geq0$ and $2\pi mx_j\leq y_j\leq2\pi nx_j$. Since $x_j=\log|w_j|$, the sequence $(x_j)$ is bounded, and so is $(y_j)$. Passing to a subsequence, we may assume that $x_j\to x\geq0$ and $y_j\to y$, where $2\pi mx\leq y\leq2\pi nx$. Hence $w=e^{x+iy}\in F_{m,n}$, proving that $F_{m,n}$ is closed. Therefore, $\Omega_{m,n}$ is open and contains the unit disc.  Set $\mu:=(\mu_m+\mu_n)/2$. Since $\mu\in C_{m,n}$ and $C_{m,n}$ is an additive semigroup, we have
\begin{equation*}
e^{s\mu}F_{m,n}\subset F_{m,n}\qquad\text{for all }s\ge0.
\end{equation*}
Consequently,
\begin{equation*}
e^{-s\mu}\Omega_{m,n}\subset\Omega_{m,n}\qquad\text{for all }s\geq0.
\end{equation*}
Thus, $\Omega_{m,n}$ is contractible and hence simply connected. Since $\mu_m\in C_{m,n}$ and $e^{\mu_m}=e$, we have $e^{-1}\Omega_{m,n}\subset\Omega_{m,n}$. Choose a Riemann map $h_{m,n}:\mathbb{D}\to\Omega_{m,n}$ satisfying $h_{m,n}(0)=0$, and define $\phi_{m,n}=h_{m,n}^{-1}\circ(e^{-1}h_{m,n})$. For $m\leq k\leq n$, we have $\mu_k\in C_{m,n}$, and hence $e^{t\mu_k}F_{m,n}\subset F_{m,n}$ for every $t\ge0$. Thus $\{m,\ldots,n\}\subset\mathscr{K}_{\phi_{m,n}}$.

Conversely, suppose that $e^{t\mu_k}F_{m,n}\subset F_{m,n}$ for every $t\geq0$. Since $1\in F_{m,n}$, for every sufficiently small $t>0$, there exist $s,r\geq0$ and $\ell\in\mathbb{Z}$ such that
\begin{equation*}
t\mu_k=s\mu_m+r\mu_n+2\pi i\ell.
\end{equation*}
Comparing real parts gives $s+r=t$. Comparing imaginary parts then yields
\begin{equation*}
kt=ms+nr+\ell.
\end{equation*}
Since $0\leq s,r\leq t$,
\begin{equation*}
|kt-ms-nr|\leq (|k|+\max\{|m|,|n|\})t.
\end{equation*}
For sufficiently small $t>0$, we must have $\ell=0$. Hence, $kt=ms+nr$. It follows that $m\leq k\leq n$. Thus $\mathscr{K}_{\phi_{m,n}}=\{m,m+1,\ldots,n\}$.
\end{example}

Unlike the preceding construction, the following example provides an explicit formula for a normalized Koenigs function that realizes any prescribed singleton set of admissible branches.

\begin{example}\label{ex:singleton}
Fix $m\in\mathbb{Z}$ and put $\mu_m:=1+2\pi i m$, $\alpha_m:=2/(1-2\pi i m)$. Consider the function
\begin{equation*}
h_m(z):=\frac{z}{(1-z)^{\alpha_m}},\hspace{1em}(z\in\mathbb{D}),
\end{equation*}
where the complex power is defined using the branch of the logarithm on the right half-plane. Clearly, $h_m(0)=0$ and $h_m'(0)=1$.

We first describe its image domain. Since
\begin{equation*}
\frac{z h_m'(z)}{h_m(z)} =\frac{1+(\alpha_m-1)z}{1-z}=1+\frac{\alpha_{m}z}{1-z},
\end{equation*}
we obtain
\begin{equation*}
\mathrm{Re}\left(\frac{1}{\mu_m}\frac{z h_m'(z)}{h_m(z)}\right)=\frac{1}{1+4\pi^2m^2}\mathrm{Re}\frac{1+z}{1-z}>0\qquad(z\in\mathbb{D}).
\end{equation*}
Thus, $h_m$ is univalent and its image domain is $\mu_m$-spirallike.

To identify the endpoint of the slit, set
\begin{equation*}
\zeta_m:=-\frac{1-2\pi i m}{1+2\pi i m}\in\partial\mathbb D,\hspace{1em}c_m:=h_m(\zeta_m)\neq0.
\end{equation*}
The point $\zeta_{m}\in\partial\mathbb{D}$ satisfies $h_{m}'(\zeta_{m})=0$. The boundary parametrization of $h_m$ shows that
\begin{equation*}
\Omega_m:=h_m(\mathbb{D})=\mathbb{C}\setminus\Gamma_m,\hspace{1em}\Gamma_m:=\{c_m e^{\mu_m s}:s\geq0\}.
\end{equation*}
Indeed, one has
\begin{equation*}
\frac{d}{d\theta}\log h_m(e^{i\theta})=\frac{\mu_m}{|\mu_{m}|^{2}}\left(2\pi m-\cot\frac{\theta}{2}\right),\hspace{1em}(0<\theta<2\pi).
\end{equation*}
Moreover, $|h_m(z)|\to\infty$ as \(z\to1\). Thus, each of the two boundary arcs joining $1$ to $\zeta_m$ is mapped onto the entire logarithmic spiral $\Gamma_m$, with endpoint $c_m$. By univalence and the boundary correspondence, the image domain is precisely $\mathbb C\setminus\Gamma_m$.

Consequently, the two boundary arcs are mapped onto the same logarithmic spiral, whose endpoint is $c_m$.

Since $e^{-t\mu_m}\Omega_m\subset\Omega_m$ for every $t\geq0$, we may define
\begin{equation*}
\phi_m(z):=h_m^{-1}\left(e^{-1}h_m(z)\right).
\end{equation*}
Then $\phi_m$ is an elliptic univalent self-map of $\mathbb{D}$ with Denjoy--Wolff point $0$ and multiplier $\phi_m'(0)=e^{-1}$; in particular, it is not an automorphism.

We claim that $\mathscr{B}_{\phi_m}^{*}=\{e^{-s\mu_m}:s\geq0\}$. Indeed, for $b\in\mathbb{C}^*$, $b\Omega_m\subset\Omega_m$ if and only if $b^{-1}\Gamma_m\subset\Gamma_m$. If this inclusion holds, then $c_m\in\Gamma_m$ implies $b^{-1}c_m\in\Gamma_m$. Hence $b=e^{-s_0\mu_m}$ for some $s_0\geq0$.

Conversely, for such a $b$,
\begin{equation*}
b^{-1}\Gamma_m=\{c_m e^{\mu_m(s+s_0)}:s\geq0\}\subset\Gamma_m.
\end{equation*}
This proves the claim.

Finally, since
\begin{equation*}
-\log\phi_m'(0)+2\pi i j=1+2\pi i j=:\mu_j,
\end{equation*}
the branch $m$ is admissible. In fact,
\begin{equation*}
\phi_t^{[m]}(z)=h_m^{-1}\bigl(e^{-t\mu_m}h_m(z)\bigr),\hspace{1em}(t\geq0)
\end{equation*}
defines a continuous one-parameter semigroup satisfying $\phi_1^{[m]}=\phi_m$.

On the other hand, suppose that $j\in \mathscr{K}_{\phi_m}$. Then $e^{-t\mu_j}\in\mathscr{B}_{\phi_m}^{*}$ for every $t\geq0$. By the characterization of $\mathscr B_{\phi_m}^{*}$, for each $t\geq0$ there exists $s_t\geq0$ such that $e^{-t\mu_j}=e^{-s_t\mu_m}$. Taking absolute values gives $s_t=t$. Therefore $e^{-2\pi i(j-m)t}=1$ for every $t\geq0$, which implies $j=m$. Hence $\mathscr{K}_{\phi_m}=\{m\}$.
\end{example}

\section{Results on the holomorphic centralizer}
Let $\phi\in\mathscr{U}(\mathbb{D})$ be an elliptic non-automorphism. Define
\begin{equation*}
 \mathscr{Z}_{\mathrm{Hol}}(\phi):=\{\psi\in\mathrm{Hol}(\mathbb{D}):\psi\circ\phi=\phi\circ\psi\}.
\end{equation*}
For a non-automorphic elliptic map, the full holomorphic centralizer admits a simple description.

\begin{prop}\label{prop:hol}
Let $\phi\in\mathscr{U}(\mathbb{D})$ be an elliptic non-automorphism with Denjoy--Wolff point $\tau$.  If $\psi\in\mathrm{Hol}(\mathbb{D})$ commutes with $\phi$, then either $\psi\equiv\tau$ or $\psi$ is univalent. In the latter case, $h\circ\psi=\psi'(\tau)h$. Consequently,
\begin{equation*}
\mathscr{Z}_{\mathrm{Hol}}(\phi)=\mathscr{Z}(\phi)\cup\{\tau\}.
\end{equation*}
\end{prop}

\begin{proof}
As before, $\psi(\tau)=\tau$. Put $a=\psi'(\tau)$ and $Q=h\circ\psi-ah$. Then $Q\circ\phi=\lambda Q$ and $Q(\tau)=Q'(\tau)=0$. If $Q$ had a zero of order $m\ge2$ at $\tau$, comparing the lowest-order non-zero terms would give $\lambda^m=\lambda$, a contradiction. Thus $h\circ\psi=ah$. If $a=0$, the univalence of $h$ yields $\psi\equiv\tau$. If $a\neq0$, then the identity $h\circ\psi=ah$ implies that $\psi$ is univalent. This conclusion also follows from [\ref{Cow'84}, Corollary 4.9].
\end{proof}

Let $\mathscr{B}_\phi:=\{a\in\mathbb{C}:a\Omega\subset\Omega\}=\mathscr{B}_{\phi}^{*}\cup\{0\}$. Theorem \ref{thm:model} and Proposition \ref{prop:hol} give the following extension.

\begin{cor}
The map $\psi\mapsto\psi'(\tau)$ is an isomorphism of topological semigroups from $\mathscr{Z}_{\mathrm{Hol}}(\phi)$ onto $\mathscr{B}_\phi$. Its inverse is given by $a\mapsto h^{-1}\circ(ah)$, where $a=0$ corresponds to the constant map $\tau$. In particular, Theorems \ref{thm:isolation} and \ref{thm:sector} remain valid with $\mathscr{Z}(\phi)$ and $\mathscr{B}_{\phi}^{*}$ replaced by $\mathscr{Z}_{\mathrm{Hol}}(\phi)$ and $\mathscr{B}_{\phi}$, respectively.
\end{cor}

\begin{remark}
The automorphic case is different.  After conjugating an elliptic automorphism to $z\mapsto\lambda z$, its holomorphic centralizer consists of the self-maps $g$ of $\mathbb{D}$ satisfying $g(\lambda z)=\lambda g(z)$.  If $\lambda$ is not a root of unity, these are exactly $g(z)=az$, $|a|\leq1$; If $\lambda$ has order $q\geq2$, then $g(z)=zF(z^q)$, where $F\in\mathrm{Hol}(\mathbb{D},\mathbb{C})$ is such that $g(\mathbb{D})\subset\mathbb{D}$. If $\lambda=1$, the automorphism is the identity map, and its holomorphic centralizer is the entire semigroup $\mathrm{Hol}(\mathbb{D})$.
\end{remark}

\section*{Acknowledgements}

The author would like to thank Professors Manuel D. Contreras, Pavel Gumenyuk, Ikkei Hotta and Ryutaro Misawa for their valuable advice during the preparation of this paper, and Professor Toshiyuki Sugawa for helpful discussions.




\end{document}